\documentclass[10pt, a4paper]{amsart}
\usepackage{amscd,amsmath,amssymb,amsthm}
\usepackage{color}
\usepackage{textcomp}

\usepackage{graphicx}

\newtheorem{theorem}{{\bfseries Theorem}}[section]

\newtheorem{proposition}[theorem]{{ Proposition}}
\newtheorem{lemma}[theorem]{{ Lemma}}
\newtheorem{corollary}[theorem]{{ Corollary}}
\newtheorem{definition}[theorem]{{ Definition}}
\newtheorem{example}[theorem]{{ Example}}
\newtheorem{remark}[theorem]{{ Remark}}

\newcommand{\bt}{\begin{theorem}}
\newcommand{\et}{\end{theorem}}
\newcommand{\bl}{\begin{lemma}}
\newcommand{\el}{\end{lemma}}
\newcommand{\bp}{\begin{proposition}}
\newcommand{\ep}{\end{proposition}}
\newcommand{\bex}{\begin{example}}
\newcommand{\eex}{\end{example}}
\newcommand{\bc}{\begin{corollary}}
\newcommand{\ec}{\end{corollary}}
\newcommand{\bo}{\begin{proof}}
\newcommand{\eo}{\end{proof}}
\newcommand{\bd}{\begin{definition}}
\newcommand{\ed}{\end{definition}}
\newcommand{\br}{\begin{remark}}
\newcommand{\er}{\end{remark}}
\newcommand{\be}{\begin{enumerate}}
\newcommand{\ee}{\end{enumerate}}

\newcommand{\T}{\mathbb{T}}

\newcommand{\cA}{{\mathcal{ A}}}

\newcommand{\cR}{{\mathcal{ R}}}
\newcommand{\cP}{{\mathcal{ P}}}
\newcommand{\cF}{{\mathcal{ F}}}

\newcommand{\cC}{{\mathcal{ C}}} 
\newcommand{\cO}{{\mathcal{ O}}}

\newcommand{\Z}{{\mathbb Z}}
\newcommand{\Q}{{\mathbb Q}}

\newcommand{\N}{{\mathbb N}}
\newcommand{\D}{{\mathbb D}}
\newcommand{\R}{{\mathbb R}}

\begin{document}

\title{Revisiting Variations in  Topological Transitivity}
\author{Anima Nagar}
\address{Department of Mathematics, Indian Institute of Technology Delhi,\\
Hauz Khas, New Delhi 110016, INDIA}

    \vspace{.2cm}
    
    \begin{abstract}
    	Topological dynamical systems $(X,T)$ are actions $T \times X \to X$, given as $(t, x) \to tx$, on a compact, Hausdorff topological space $X$ with $T$ as an acting group or monoid. We take up the property of topological transitivity especially for semiflows $(X,S)$ and discuss the variations in its definitions. We emphasize on those properties  of transitivity that differ for semiflows  as compared to those for flows or  (semi)cascades. \end{abstract}

\maketitle

\renewcommand{\thefootnote}{}

\footnote{\emph{keywords:} minimality, topological transitivity, strong transitivity, very strong transitivity, locally eventually onto, strongly product trantivity, semiflows}


\footnote{ The author thanks the   anonymous referees for very detailed comments  and suggestions, that has tremendously polished the entire write-up.}

\medskip

 	``\emph{Transitivity}'' gets its \emph{form} from the concept of ``\emph{recurrence}''. The frequencies of such recurrence in transitive systems were first explored and utilized by \textbf{Hillel Furstenberg}  especially in the field of combinatorial number theory, where he used \emph{`families'} to characterize recurrence.

 	In his seminal paper on disjointness in topological dynamics, \textbf{Furstenberg} \cite{f}  laid a foundation for the classification of dynamical systems by their recurrent properties. He systematized the theory of recurrences by considering (Furstenberg) families, in particular defining such families for classifying transitivities. We follow this tradition and look into more details on various transitivities that were discussed in \cite{aan}.

 Based on the study developed by Furstenberg, \textbf{Akin} \cite{ak} extended the notion of recurrence  to the setting of general monoid actions. In this context, a `\emph{Fursterberg family}' on the subsets of an acting monoid, is considered to define the recurrence of its action.

 \bigskip

In this article, we discuss some variations of topological transitivity for the systems $(X, T)$, on the lines of the study made in \cite{aan}, infused by the ideas in \cite{ak,f} with emphasis on the properties of ``strongly transitive'' systems.

For a more general approach, the survey \cite{ckn} is recommended.

\begin{section}{ Topological Dynamical Systems}

A pair $(X,T)$ is called a \emph{topological dynamical system}, where
 $X$ is a compact, Hausdorff topological space and $T$ is a topological group or monoid or semigroup acting on $X$.  We say that $T$ acts on $X$ when there is a continuous map $\pi: T \times X \to X$ such that

\be

\item $\pi (e, x) = x$ for all $x\in X$, if $e \in T$ is the identity in the monoid or group $T$;

\item $\pi (ts, x) = \pi (t, \pi (s,x))$ for all $t, s \in T$ and
all $x\in X$.

\ee

\medskip

 $X$ is called the \emph{phase space}, $T$ the \emph{acting group or monoid} and the \emph{action} $\pi$ gives the homeomorphism or continuous map $\pi^t: X \to X$ defined as $\pi^t(x) = tx$.

\br When $T$ is a topological group - $(X,T)$ is called a \emph{flow}, and for a topological monoid $S$ - we call $(X,S)$ a \emph{semiflow}. When the acting group $T = \Z$, then $\pi^1 = f$ gives a generating homeomorphism on $X$, i.e. $f(x) = \pi(1,x)$, giving iterations $f^n(x) = \pi(n,x) = nx$. In this case we call the system $(X,f)$ a \emph{cascade}. A \emph{semicascade}  $(X,f)$ where $f: X \to X$ is a continuous mapping, corresponds to the case when the monoid $S = \Z_+ = \N \cup \{0\}$.

In general, we call $(X,T)$ - a system without specifying if it is a flow, semiflow, cascade or semicascade. We  shall specify  specific cases, though our results are mostly for monoids.

Also our results remain true for  any particular topology on $T$, unless we specify more assumptions on this topology.
\er

\medskip

The \emph{orbit} of $x \in X$ under $T$ is $Tx$, which is the smallest \emph{$T$-invariant} set containing $x$. For a cascade $(X,f)$ we also call   $\cO(x) = \{f^n(x) : n \in \Z_+\}$  the \emph{(forward)orbit} of the point $x$. We call $ {\overline{Tx}} $ (or $\overline{\cO(x)}$ for cascade)  the \emph{orbit closure} of $x$ under $T$ and it
is the smallest closed $T$-invariant set containing $x$.  A set
$A \subset X$ is called \emph{$T$-invariant} if $TA\subset A$, where $TA = \{ta: t \in T, \ a \in A\}$.

\medskip

 For a  cascade or semicascade $(X,f)$, $x_0 \in X$ is called a \emph{fixed point} if $f(x_0) = x_0$. And $y_0 \in X$ is called a \emph{periodic point} if there exists $n \in \N$ such that $f^n(y_0) = y_0$. The smallest such $n$ is called the \emph{period} of $y_0$.

\medskip

 Let $(X,T)$ and $(Y,T)$ be dynamical systems.  A continuous, surjective map
$\phi \colon X \to Y$ is called a \emph{factor-map or semi-conjugacy} if $\phi(tx) = t\phi(x)$ for all $x\in X$ and $t\in T$. In addition if $\phi$ is a homeomorphism, we
say that  $(X,T)$ and $(Y,T)$ are \emph{conjugate} as dynamical systems.

\medskip

For a (semi)cascade, the \emph{$\omega$-limit set} of a point $x \in X$ under $f$ is the set of all limit points of $\{f^n(x): n \in \Z (\Z_+)\}$. Thus $y \in \omega(x)$ if and only if there exists sequence $\{n_k\} \nearrow \infty$ such that $f^{n_k}(x) \to y$ or a net $\{n_k\}$ on some increasing directed set.

 The \emph{$\omega$-limit set} of a point $x \in X$ in a (semi)flow $(X,T)$, denoted as
$\omega(x)$, is the set of all limit points of
$\{tx: t \in T\}$. Thus we have that for $y \in \omega(x)$ and any non-empty, open $U \ni y$, the set $\{t: tx \in U\}$ is unbounded in $T$ or has no compact closure in $T$.

\br Observe that for all $t\in T$ and $x \in X$, we have $t\omega(x)\subset \omega(x) \subset \omega(tx)$. \er

 For very nice (say sigma-compact and locally compact) $T$, we can say that $y \in \omega(x)$ if  for  every  open  $U$ with $y \in U$, there exists compact $K_1 \subset K_2 \subset \ldots \subset K_i \subset \ldots \subset T$  such that $s_ix \in U$ for $s_i \in T \setminus  K_i$, i.e. $x \in  \bigcap \limits_{s_i \in T \setminus  K_i}  s_i^{-1}U$.

\medskip

We note that here we have just taken a particular representation of an ``admissible set''. In general such admissible sets could vary. Some of such admissible sets were first considered by \textbf{Furstenberg} \cite{f}.

\medskip

We denote by $\cP = \cP(T)$ the set
of all subsets of $T$. A subset $\cF \subset \cP$  is a \emph{Furstenberg family}, if  it is hereditary
upwards, that is, $F_1 \subset F_2$ and $F_1 \in \cF$ implies $F_2 \in \cF$. The family $\cF$ is \emph{ proper} if it is a proper subset of $\cP$. Any $\cA \subset \cP$ clearly generates a family $\{ F \in \cP : F \supset A$ for some $A \in \cA \}$.
A \emph{filter} $\cF$ is a proper family closed under intersection, that is, $\cF$ is a proper subset of $\cP$ and for $F_1, F_2 \in \cF$ implies $F_1 \cap F_2 \in \cF$. The family $\cF$ is a filter if and only if $\cF$ is a proper family and $\cF \cdot \cF \subset \cF$. A family $\cF$ has the\emph{ Ramsey property} if $F \in \cF$ and $F = F_1 \cup F_2$ imply that $F_i \in \cF$ for some $i \in \{1, 2\}$.

For a family $\cF$, the \emph{dual family of $\cF$}, denoted by $ \cF^*$, is defined as
$\cF^* = \{F \in \cP : F \cap F' \neq \emptyset,$ for any $F' \in \cF\}$.

\bigskip

For families $\cF_1, \cF_2 \subset \cP$, we define $\cF_1 \cdot \cF_2 = \{ F_1 \cap F_2: F_1 \in \cF_1, \ F_2 \in \cF_2\}$.

\medskip

We call a family $\cF$ to be \emph{translation invariant} if  $ F \in \cF$ if and only if $t^{-1}(F) \in \cF$ (where $t^{-1}(F)$
denotes the preimage of $F$ by $t$), for all $ t \in T$. Also $\cF$ is \emph{thick} if for $F \in \cF$ and $t_1, \ldots, t_k \in T \Leftrightarrow \bigcap \limits^k_{i=1} \ t_i^{-1}(F) \in \cF$.

\medskip

Now $F \subset T$ is called \emph{syndetic} if there is a compact $K \subset T$ such that $T = KF =\{kf: k \in K $ and $ f \in F \}$. A set $F \subset \Z (\N)$ is called syndetic if it is relatively dense, i.e., it does not contain arbitrarily large gaps.

We call a family $\cF$ to be \emph{syndetic} if every  $ F \in \cF$ is syndetic.

\medskip

A set $F_1$ is syndetic if it intersects each thick set and a set $F_2$ is thick if it intersects each syndetic set.

\medskip 

Also $F \subset T$ is called \emph{piecewise syndetic} if it is an intersection of  thick and syndetic subsets of $T$. We call a family $\cF$ to be \emph{piecewise syndetic} if every  $ F \in \cF$ is piecewise syndetic. Thus for syndetic family $\cF_1 \subset \cP$ and thick family $\cF_2 \subset \cP$, the family $\cF = \cF_1 \cdot \cF_2$ is piecewise syndetic.

\medskip

Let $\star$ be the binary operation for the group or monoid $T$.  Define the finite sums of a sequence $\{t_i\}$ in $T$ as

$$FS\{t_i\}=\left \{ t_{i_1} \star t_{i_2} \star \ldots \star t_{i_n}: n \in \N \right\}.$$

The set $F \subset T$ is called an \emph{IP set} if there exists a sequence $\{t_i\}$ in $T$ such that $FS\{t_i\}\subset F$. A family $\cF$ is called an \emph{IP family} if $\cF$ is the family of all IP sets.

\medskip

We recommend the excellent treatment by \textbf{Akin} \cite{ak} for a detailed account on these admissible families.

\medskip

 \textbf{Furstenberg} defined families and classified different types of transitivities based on the combinatorial properties of the families of hitting sets \cite{f}.

\medskip

We identify  a singleton with the point it contains.  For any two non-empty, open $U,V \subset X$ and $x \in X$, for the system $(X,T)$, where $T$ is any group or monoid, we define the hitting times:
$$ N_T(x,V) \quad =  \ \{ t \in T : t(x) \in V  \} = \ \{ t \in T : x \in t^{-1}(V)  \}.$$
$$ N_T(U,V) \quad =  \ \{ t \in T : t(U) \cap V \not= \emptyset \} = \ \{ t \in T : t^{-1}(V) \cap U \not= \emptyset \}.$$
$$ N_T(U,x) \quad = \quad \{ t \in T : x \in t(U) \} \ = \ \{ t \in T : t^{-1}(x) \cap U \not= \emptyset \}.$$

We see that for $x \in U$, $y \in V$ we have $N_T(x,V) \subset N_T(U,V)$ and $N_T(U,y) \subset N_T(U,V)$.

\medskip

For cascades or semicascades, we suppress the subscript $T$ from the hitting times and denote these by $N(x,V)$, $N(U,V)$ and $N(U,x)$ and consider only positive instances as first defined in \cite{f}.

\medskip

Let $S$ be any semigroup or monoid and $[S]$ denote the smallest group generated by $S$. Let $S^{-1} \subset [S]$ denote the set of inverses of elements of $S$.

We note that $ N_S(x,V) = \ \{ s \in S : s(x) \in V  \} = \ \{ s \in S : x \in s^{-1}(V)  \} = \ \{ s \in S^{-1} : x \in s(V) \} \ = N_{S^{-1}}(V,x)$. Our hitting time sets are interrelated.

\medskip

We consider the product (semi)flow $(X^n,T)$ as $(x_1,x_2, \ldots, x_n) \stackrel{t}{\rightarrow} (tx_1,tx_2, \ldots, tx_n)$ for any $n \in \N$.

\bigskip

For  each $F \in \cP$, every point $x \in X$ and each non-empty, open  $U \subset X$  define the \emph{$F-$orbit} as $T^{F} x = \{tx : t \in F\}$. The $\omega$-limit set of $x$ with respect to $\cF$ ,
denoted by $\omega_{\cF}(x)$, is defined as $\omega_{\cF}(x) = \{y \in X : N_T (x,W) \in  \cF$ for every open $W \ni  y\}$.

\medskip

\medskip

\bl \label{iphit} Let $T$ be any Abelian group or monoid. If $(y, y) \in  \omega((x,y))$, then for every non-empty, open $U \ni y$, the hitting time set $N_T(x, U)$ is an IP set. \el

\bo We prove this by induction. For non-empty, open $U \ni y$, let $U_1 = U$ then there exists $t_1 \in T$ such that $t_1 (x,y) \in U_1 \times U_1$. Let $U_2 = U \cap t_1^{-1}U_1$, then since $y \in U_2$  there exists $t_2 \in T$ such that
$t_2 (x,y) \in U_2 \times U_2$. 

Then for every $m \in FS\{t_1, t_2\}$, we have $mx \in U$.

Suppose that we get a sequence $\{t_1, t_2, \ldots, t_{n-1}\}$ in $T$ and a sequence of non-empty, open sets $U = U_1 \supset U_2, \ldots , U_{n-1}$ , such that for every $m \in FS\{t_1, t_2, \ldots t_{n-1}\}$, we have $m(x,y) \in U \times U$.

Let $U_n = U \cap \left(\bigcap \limits_{m \in FS\{t_1, t_2, \ldots , t_{n-1}\}} m^{-1}(U)\right)$.

Then $U_n \ni y$ and so there exists $t_n \in T$ with
$t_n(x,y) \in U_n \times U_n$.

This gives for every $m \in FS\{t_1, t_2, \ldots, t_{n}\}$, we have $m(x,y) \in U \times U$.

Thus, by induction we get a sequence $\{t_i\} \subset T$ for which $FS\{t_i\} \subset N(x, U)$.
\eo

\medskip

   A point $x \in X$ is said to be
\emph{non-wandering} in $(X,T)$ if for every neighbourhood $U$ of $x$, there is a $t
\in T$ such that $t(U) \cap U \neq \emptyset$. The set of all
non-wandering points of $(X,T)$ is denoted as $\Omega(T)$.

For a cascade or semicascade, a point $x \in X$ is
\emph{non-wandering} if for every neighbourhood $U$ of $x$ there is a $n
\in \mathbb{N}$ such that $f^n(U) \cap U \neq \emptyset$. The set of all
non-wandering points of $f$ is denoted as $\Omega(f)$.

\br Note that $x$ is a non-wandering point if for every open  $U \subset X$   with $x \in U$,  $N_T(U,U) \neq \emptyset$. \er

For any Furstenberg family $\cF$, a   point $x \in X$ is called an \emph{$\cF$-non-wandering point} if  for every open  $U \subset X$   with $x \in U$,  one has $N_T(U,U) \in \cF$. The set of all $\cF$-non-wandering points, denoted by $\Omega_{\cF}(T)$, is defined as $\Omega_{\cF}(T) = \{x \in X : N_T (U,U) \in  \cF$ for every open $U \ni  x\}$.

\medskip

 A point $x \in X$ is called \emph{recurrent} whenever $x \in \omega (x)$. The set of all
recurrent points of $(X,T)$ is denoted as $\cR(X)$.

\br Also note that $\cR(X)$ is $T$-invariant, i.e., $t(\cR(X)) \subset \cR(X)$ for all $t \in T$. \er

For any Furstenberg family $\cF$, a   point $x \in X$ is called \emph{$\cF$-recurrent} if and only if  $x \in \omega_{\cF}(x)$.

\medskip

\bd   A point $x \in X$ is said to be
\emph{strongly non-wandering} in $(X,T)$ if for every neighbourhood $U$ of $x$, there is a $t
\in T$ such that $x \in t(U) $. The set of all
strongly non-wandering points of $(X,T)$ is denoted as $\widehat{\Omega}(T)$.

For a cascade or semicascade, a point $x \in X$ is
\emph{strongly non-wandering} if for every neighbourhood $U$ of $x$ there is a $n
\in \mathbb{N}$ such that $x \in f^n(U)$. The set of all
strongly non-wandering points of $f$ is denoted as $\widehat{\Omega}(f)$. \ed

\br Note that $x$ is a strongly non-wandering point if  for every open  $U \subset X$   with $x \in U$, one has $N_T(U,x) \neq \emptyset$. \er

\br Note that $x$ is a strongly non-wandering point, then it is also non-wandering, i.e., $\widehat{\Omega}(T) \subset {\Omega}(T)$. \er

\bd For any Furstenberg family $\cF$, a   point $x \in X$ is called an \emph{$\cF$-strongly non-wandering point} if  for every open  $U \subset X$   with $x \in U$, one has  $N_T(U,x) \in \cF$. The set of all $\cF$-strongly non-wandering points, denoted by $\widehat{\Omega_{\cF}}(T)$, is defined as $\widehat{\Omega_{\cF}}(T) = \{x \in X : N_T (U,x) \in  \cF$ for every open $U \ni  x\}$.\ed

\medskip

\br Note that for systems $(X,T)$ and $(Y,S)$ and Furstenberg families $\cF_T \subset \cP(T)$ and $\cF_S \subset \cP(S)$, we have

$${\Omega_{\cF_T}}(T) \times {\Omega_{\cF_S}}(S) \subset {\Omega_{\cF_T \cdot \cF_S}}(T \times S),$$

$$\widehat{\Omega_{\cF_T}}(T) \times \widehat{\Omega_{\cF_S}}(S) \subset \widehat{\Omega_{\cF_T \cdot \cF_S}}(T \times S),$$

\smallskip

and for systems $(X,T)$ and $(Y,T)$ and Furstenberg family $\cF = \cF_X = \cF_Y \subset \cP(T)$ considering actions on X and Y resp., we have

$${\Omega_{\cF_X}}(T) \times {\Omega_{\cF_Y}}(T) \subset {\Omega_{\cF_X \cdot \cF_{Y}}} (T \times T),$$

$$\widehat{\Omega_{\cF_X}}(T) \times \widehat{\Omega_{\cF_Y}}(T) \subset \widehat{\Omega_{\cF_X \cdot \cF_{Y}}} (T \times T),$$

\smallskip

with equality when $\cF$ is a filter.
\er

\medskip

We call the system $(X,T)$ to be \emph{central} if each $t \in T$ is surjective, i.e., $t(X) = X$ for all $t \in T$.

\medskip

The system $(X,T)$ is called \emph{$\cF$-central}  for a family $\cF$ if  $N_T(U, U) \in \cF$ for all $x  \in X$ and open $U \ni x$. This essentially means that $\Omega_{\cF}(T) = X$.

\medskip

The system $(X,T)$ is called \emph{$\cF$- strongly central}  for a family $\cF$ if  $N_T(U, x) \in \cF$ for all $x  \in X$ and open $U \ni x$. This essentially means that $\widehat{\Omega_{\cF}}(T) = X$.

\medskip

A subset $S$ in an Abelian $T$ is said to be \emph{replete} if $S$ contains some translate of each
compact set in $T$. A set $A \subset T$ is said to be \emph{extensive} if $A$ intersects
every replete semigroup in $T$. A point $x \in X$ is said to be recurrent under $T$ if for every neighborhood $U \ni x$
 there corresponds an extensive set $A$ in $T$ such that $Ax \subset U$. This study of recurrence was first taken up by  \textbf{Gotschalk and Hedlund} \cite{goh}. For non-Abelian groups the existence of replete semigroups seems to be quite rare, and so the definition of recurrence in these terms appears to be inadequate. But the basic idea in recurrence still keeps this structure of some admissible set of repeats of occurrences for the point that is recurrent.

\medskip

We refer to \cite{ak, f, join, gh, v1} for more details on transformation systems.

\end{section}

\begin{section}{Various forms of Topological Transitivity}

\begin{subsection} {Minimal Systems}

Let $(X,T)$ be a flow. The simplest dynamics that one can observe is when the system is ``minimal''. The best treatment of minimal flows is by \textbf{Auslander} \cite{aus}.

\medskip

A set $M \subset X$ is called a \emph{minimal set} if $M$ is closed, non-empty and $T-$invariant and $M$ has no proper subset with these properties, i.e. if $N \subseteq M$ is closed and invariant then $N = M$ or $N = \emptyset$.

\br Note that $M \subset X$ is minimal if and only if it is the orbit closure of each of its points, i.e., $\forall x \in M$, we have $M = {\overline{Tx}}$. \er

If $X = {\overline{Tx}}$ $\ \forall x \in X$, then the flow $(X,T)$ is called a \emph{minimal flow}.

\medskip

A point $x \in X$ is called an \emph{almost periodic point} if for every neighbourhood $U$ of $x$, there is a syndetic  $A \subset T$ such that $Ax \subset U$.

\br Note that $x$ is an almost periodic point if for every  open  $U \subset X$ with $x \in U$, one has $N_T(x,U) $ is syndetic. \er

We refer to \cite{aan} for more details on minimal cascades.

\medskip

Noninvertible minimal semicascades were studied by \textbf{Kolyada,  Snoha  and  Trofimchuk} in \cite{kst}.

\bex  Let $\Lambda= \{0,1\}$  and define $X=\Lambda ^\N$. We consider the shift map $\sigma: X \to X$.

To obtain a minimal subset of $X$, it is enough to construct an almost periodic point $p \in X$ since then ${\overline{\cO(p)}}$ will be minimal. We mimic the classical construction due to Marston Morse and   Axel Thue, giving the Morse-Thue sequence.

This construction is done using substitution: $0 \to 01, \ 1 \to 10$. Hence,
$$0 \to 01 \to 0110 \to 01101001 \to 0110100110010110 \to \cdots$$
This will finally converge to some $x \in \{0,1\}^\N$. This construction indicates that every finite word in $x$ occurs syndetically often.. Then $x$ is almost periodic in $ \{0,1\}^\N$, and so $({\overline{\cO(x)}}, \sigma)$ is a noninvertible minimal dynamical system. \eex

A map $f: X \to X$ is called \emph{almost open} if for every non-empty open subset $U \subset X$, the set $f(U)$ has non-empty interior in $X$. Essentially a map is almost open if and only if the inverse
image of every dense subset is dense. A map $f : X \to Y$ is called \emph{irreducible} if the only closed  $A \subset X$ for which $f(A) = Y$ is $A = X$. In particular, an  irreducible map is always  surjective. Also, irreducible maps are always almost open. Minimal maps are irreducible, and irreducible  maps are almost open; thus implying minimal semicascades to be almost open.

\medskip

\br We note that irreducibility or almost openess does not imply minimality even for cascades as can be seen in the below example.

Let $X = \T \times \{1,2\}$ where $T$ is  the unit circle bijective with $[0,1)$ and consider the irrational rotation on both circles in $X$ given as $\tau (\theta, i) = (\alpha \theta, i)$ for $i=1,2$. Then $(X, \tau)$ is irreducible, almost open but not minimal.\er

For more on properties of minimal semicascades we refer to \cite{kst}.

%
%
%
%
%
%
%

\medskip

\textbf{This leads to the natural question on the properties of minimal semiflows.}

We note that recently such a study is done by \textbf{Auslander and Dai} \cite{ad}. However, they have a different perspective.

\bd A \emph{minimal semiflow} is the system $(X,S)$, where $S$ is a semigroup or monoid action on $X$, for which  $X = {\overline{Sx}}$ is true $\ \forall x \in X$. \ed

\medskip

\bex \label{nao} 
	Consider  $X = \T$, the unit circle bijective with $[0,1)$ and consider
	
	$$S = \{f_{p,q}(x) := px+q \ (\mod 1): p \in \Q \cap (0,1], q \in \Q \cap [0,1]\}$$
	
	with operation as composition of functions. Then $S$ is an equicontinuous family in $\mathcal{C}(X)$ (the space of all continuous real valued functions on $X$ with the uniform topology), and each  $ f_{p,0} \in S $ is invertible but   $ f_{p,0}^{-1} \notin S $ for all $ p \in (0,1) $.  Also the monoid $ S $ acts on $ X $ and $ (X,S) $ is a semiflow.
	
	\medskip
	
	For the semiflow $(X,S)$, we note that $\overline{Sx} \supset \overline{\{f_{1,q}(x): q \in [0,1] \cap \Q\}}=  X$ for all $x \in X$ and so $(X,S)$ is minimal.
	
	\medskip

	Let $0 \in X$ and open $U \ni 0$ and let $A = \{f_{r,s} \in S: f_{r,s}(0) \in U\} = N_S(0,U)$.

	Note that if $A$ is syndetic then there would exist compact $K \subset S$ such that  $S = KA$. 
	
	Observe that if  $U = [0,\frac{1}{27}) \cup (\frac{26}{27},1]$ with $0$ identified with $1$, then here
	
	$$A =  \left\{f_{r,s}:  \ r \in (0,1] \cap \Q, \ s \in \left( [0,\frac{1}{27})  \cup  (\frac{26}{27},1]  \right) \cap \Q \right\}. $$

	Let $h=\frac{1}{81}$ then	$ f_{h,0}S \subset A $. But $ f_{h,0}^{-1} \notin S $ and so we cannot have $ S = f_{h,0}^{-1}A $. 
	
	If there is a compact $K \subset S$ such that $S=KA$ then for all $s \in \Q \cap [0,1]$, there exists $f_{1,k} \in K$ and $f_{1,a} \in A$ such that $f_{1,s} = f_{1,k}f_{1,a}$, i.e. $s=k+a$ for all $s \in \Q \cap [0,1]$. Since $f_{1,a} \in A$ only for $a \in  \left( [0,\frac{1}{27})  \cup  (\frac{26}{27},1]  \right) \cap \Q$,  we must have $f_{1,t} \in K$, say  for all $t \in [\frac{5}{27},\frac{13}{27}] \cap \Q$. Let $i \in [\frac{5}{27},\frac{13}{27}] \cap (\R \setminus \Q)$and let $r_n \in [\frac{5}{27},\frac{13}{27}] \cap \Q$ be such that $r_n \to i$. If $K$ is compact, $\{f_{1,r_n}\}$ must uniformly converge in $K$. This is not possible. 
	
	\medskip

	Thus $N_S(0,U) $ is not syndetic.

\eex

\br For a minimal semiflow $(X,S)$, each $s \in S$ need not be surjective, nor almost one-to-one, nor  $N_S(x,U) $ be syndetic for every  open  $U \subset X$ with $x \in U$.

However, Zorn's lemma guarantees that every semiflow will have a minimal closed invariant subset. \er

\begin{remark}
	We note that in the example above, every $s \in S$ is not surjective. Hence though every point of $X$ is recurrent, the system is not distal. In fact the semiflow $(X,S)$ is proximal. 
\end{remark}

We note here that if  $[S]$ denote the smallest group generated by the monoid $S$, then $ (X, S)$ proximal implies that $ (X,[S]) $ is also proximal. However the converse here may not be true. For example we can consider $S^{-1}$ in Example \ref{nao}. Note that  $[S] = \{f_{a,b} := ax+b: a \in (0,\infty) \cap \Q, \ b \in \Q\}$ here.

\bigskip

\br It has been observed in \cite{sg} that the only endomorphism a minimal proximal flow admits is the identity automorphism. We note that minimal proximal semiflows also admit  the identity automorphism as the only endomorphism, and the proof is similar as the one for flows.

However, every minimal proximal
flow is weakly mixing \cite{sg} whereas Example \ref{nao} is a minimal proximal semiflow that is not weakly mixing. In fact, it is equicontinuous. \er

We summarize the trivial properties of minimal sets in semiflows:

\bp For  the semiflow $(X,S)$:

\be

\item If $M_1$ and $M_2$ are minimal subsets of $X$ for any semiflow $(X,S)$ then either $M_1 = M_2$ or $M_1 \cap M_2 = \emptyset$.

\item Let  $(X,S)$ be a semiflow. Then $X$ contains a minimal set.

\item If $(X,S)$ is minimal then the only closed, invariant subsets of $X$  are $\emptyset$ and $X$.

\item For semiflows $(X,S)$ and $(Y,S)$, let $\pi: X \to Y$ be a semi-conjugacy.

Then if $X_0 \subset X$ is minimal then $\pi(X_0) = Y_0 \subset Y$ is minimal.
\ee
\ep

\medskip

\bd Consider a semiflow $(X,S)$. The semiflow $(X,S)$ is \emph{irreducible} if  every $s \in S$ is irreducible (according to the definition in \cite{kst}). \ed

We note that the image of no proper closed subset of $X$ under the action of the semigroup $S$ can be equal to $X$ in an irreducible semiflow. Also when  $(X,S)$ is irreducible then it is also central.

\br We note that in the Example \ref{nao}, the semiflow $(X,S)$ is not irreducible. Thus for semiflows, minimality does not imply irreducibility. \er

\bd \label{open} Consider a semiflow $(X,S)$. Then $s \in S$ is almost open if $sU$  has a non-empty interior for all non-empty, open $U \subset X$. Equivalently $s^{-1}(D)$is dense in $X$ whenever $D$ is dense in $X$.  The semiflow $(X,S)$ is \emph{almost open} if each $s\in S$ is almost open.

The semiflow $(X,S)$ is \emph{ open} if each $s\in S$ is  an open map.\ed

\br We note that in  Example \ref{nao}, all the elements in $S$ are almost open. Thus there are semiflows that are not irreducibe but almost open. \er

\bp \label{mini} Let a semiflow $(X,S)$ be minimal and $U \subset X$ be non-empty and open. Then there exists $s_1, \ldots, s_m \in S$ such that $X = \bigcup \limits_{j=1}^m {s_j}^{-1}(U) $.

Furthermore, if $S$ is Abelian and  also central, then there exists $s_{n_1}, \ldots, s_{n_m} \in S$ such that $X = \bigcup \limits_{j=1}^m {s_{n_j}}(U) $.
\ep

\bo Since $(X,S)$ is minimal, $\overline{Sx} = X$ for all $x \in X$.  Consider a non-empty, open set $U \subset X$. For every $x \in X$, there exists $s \in S$ such that  $s(x) \in U$. And so we have open neighbourhood $V_x \ni x$ with $ V_x \subset s^{-1}(U)$. Then the class of open sets $\{ V_x : \  \ x \in X\}$  forms an open cover of $X$. By compactness of $X$, this open cover has  a finite sub-cover and so $X = \bigcup \limits_{j=1}^m V_{x_j}$. Also corresponding to each $x_j$ there is an $s_j \in S$ such that $V_{x_j} \subset s_j^{-1}(U)$. And so $X =  \bigcup \limits_{j=1}^m V_{x_j} \subset \bigcup \limits_{j=1}^m {s_j}^{-1}(U) \subset X$.

\medskip

Furthermore when $S$ be Abelian with each $s \in S$ surjective, then $s(X) = X$, $\forall \ s \in S$. Now
$X =   \bigcup \limits_{j=1}^m {s_j}^{-1}(U) $. Then for $s = s_1 \circ \ldots \circ s_m$ if  $s_{n_j} = s{s_j}^{-1}$,

 $X = s(X) =   \bigcup \limits_{j=1}^m s{s_j}^{-1}(U) = \bigcup \limits_{j=1}^m {s_{n_j}}(U). $\eo

\begin{corollary} Let a semiflow $(X,S)$ be minimal and irreducible with $S$ Abelian. Then for all non-empty, open $U \subset X$  then there exists $s_{n_1}, \ldots, s_{n_m} \in S$ such that $X = \bigcup \limits_{j=1}^m {s_{n_j}}(U) $.
\end{corollary}

\bigskip

\end{subsection}

\begin{subsection}{Topological Transitivity }

Topological transitivity can be described as the eventuality of  neighbourhood of every point to visit every
region of the phase space at some time.

 A flow or semiflow $(X, T)$ is called \emph{topologically transitive} or \emph{transitive} if for all non-empty open sets $U, V \subset  X$, there is a  $t \in T $ for which $t(U) \cap V \neq \emptyset$ [Equivalently,  $U \cap t^{-1}(V)  \neq \emptyset$].

 A cascade or semicascade $(X,f)$ is
said to be  \emph{topologically transitive} if for every pair of non-empty
open sets $U,V$ in $X$, there is a $n \in \mathbb{N}$
such that $f^n(U) \cap V \neq \emptyset$ [Equivalently,  $U \cap f^{-n}(V)  \neq \emptyset$].

\medskip

A flow or semiflow $(X, T)$ is called \emph{point transitive} if
there is an $x_0 \in X$ such that $\overline{Tx_0} = X$ ( i.e.
$\{tx_0: t \in T\}$ is  dense in $X$).

An $x \in X$ is called a \emph{transitive point}, if $\overline{Tx} = X$. We note that $x$ is transitive when $\omega(x) \cup Tx$ = X.

We denote the set of all transitive points of $(X,T)$ by $Trans(T)$.

\br We note that for a minimal system $Trans(T) = X.$ \er

The cascade or semicascade $(X,f)$ is said to be \emph{point transitive} if there is an $x_0 \in X$ such that $\overline{\cO(x_0)} = X$ ( i.e.
$X$ has a dense orbit). All such points with dense orbits are called \emph{transitive points} and the set of transitive points in $X$ is denoted as $\emph{Trans(f)}$.

\medskip

Both these definitions of point transitivity and topological transitivity for (semi)cascades are equivalent, in a wide
class of spaces, including all perfect, compact metric spaces. We recall,

\bp  \cite{k} For a cascade or semicascade $(X,f)$, if $X$ has no isolated point then  point transivity  implies the transitivity of $(X,f)$.

The converse holds if $X$ is separable and of second category. \ep

Though we see that this is not essentially true for a semiflow $(X,S)$ when $S$ is not discrete.

\bex \label{sgex} Let $S = (\R^+, +)$ be the monoid of nonnegative real numbers  and $X = [0, \infty]$, the
 one-point compactification of the nonnegative reals.
Define $S \times X \to X$, by $(s, x) \to s + x$. Then for the natural action of $S$ on $X$, $(X,S)$ is a semiflow.

  Then for the semiflow
$(X,S)$, $\overline{S0}$ is dense in $X$ and so $Trans(S) \neq \emptyset$. Thus $(X,S)$ is point transitive but it can be clearly seen that $(X,S)$ is not topologically transitive.
\eex

\medskip

\br We note this obvious fact: if  $(X,S)$ is minimal, then $(X,S)$ is topologically transitive. \er

\medskip

The proof of the following follows trivially from known facts.

\bp We have the following   equivalent conditions for transitivity of semiflows $(X,S)$:

\be
\item $(X,S)$ is topologically transitive.

\item for every pair of non-empty open sets $U$ and $V$ in $X$, there
is a  $s \in S$ such  that $s^{-1}(U) \cap V \neq
\emptyset$.

\item for every pair of non-empty open sets $U$ and $V$ in $X$, $N_S(U,V) \not= \emptyset$.

\item for every non-empty open set $U \subset X$, the set
${\cup}_{s \in S} s(U)$ is dense in $X$.

\item for every non-empty open set $U \subset X$, the set
${\cup}_{s \in S} s^{-1}(U)$ is dense in $X$.

\item if $E \subset X$ is closed and $SE \subset E$, then either $E = X$
or $E$ is nowhere dense in $X$.

\item if $U \subset X$ is open and $U \subset SU$, then either $U =
\emptyset$ or $U$ is dense in $X$.

\ee

\medskip

Moreover, if X is a compact, perfect metric space, then

\be

\item There exists $x \in X$ such that the orbit $Sx$ is dense in $X$,
i. e. the set $Trans(S)$ of transitive points is non-empty.

\item The set $Trans(S)$ of transitive points equals $ \{ x : \omega (x) = X \}$ and it
is a dense  subset of $X$.
\ee
\ep

\medskip

\bex \label{tm} Let $f:[0,1]
\rightarrow [0,1]$ be the \emph{tent map} defined as $f(x) = \left\{
                                           \begin{array}{ll}
                                             2x, & \hbox{$0 \leq x \leq 1/2$;} \\
                                             2(1-x), & \hbox{$1/2 \leq x \leq 1$.}
                                           \end{array}
                                         \right.$ Then the semicascade
$([0,1],f)$ is transitive. 

Here for any non-empty open $J$ in $[0,1]$, $ \exists \ n \in \mathbb{N}$ such
that $f^n(J) = [0,1]$. Thus $([0,1], f)$ is transitive.

 \eex

\br Any transitive homeomorphism is irreducible, almost open and almost one-to-one. \er

For more on transitivity, we refer to \cite{aan, gh, k, nk}.

\medskip

For a Furstenberg family in $\N$, \emph{$\cF-$point transitivity} was first defined by  \textbf{Jian Li}  in \cite{j} for (semi)cascades. In general, $\cF-$transitivity does not imply $\cF-$point transitivity.
\medskip

We look into this concept for (semi)flows. Let $\cF$ be a Furstenberg family of subsets of $T$. A point $x \in X$ is called
an \emph{$\cF-$transitive point} if for every non-empty, open $ U \subset X $, one has $N_T(x,U) \in \cF$.

The (semi)flow $(X, T )$   is called \emph{$\cF-$point transitive} if
there exists some $\cF-$transitive point $x_0 \in X$.

\medskip

For the Furstenberg family $\cF$, the (semi)flow  $(X,T)$ is called \emph{$\cF$-transitive} if for every
pair $U, V$ of non-empty, open sets in $X$, one has $N_T(U,V) \in \cF$.

\br Note that if $(X,T)$ is  \emph{$\cF$-transitive} then $(X,T)$ is also \emph{$\cF$-central} and $\Omega_\cF(T)= X$ though it need not be $\cF-$point transitive. \er

We recommend the enthusiastic reader to look into a detailed description of $\cF$-transitivity discussed by \textbf{Akin} in \cite{ak}.

\bigskip

\end{subsection}

\begin{subsection}{Weakly Mixing and (Strongly) Mixing}

 A (semi)flow  $(X,T)$ is said to be \emph{weakly mixing} if the product (semi)flow  $(X^2, T)$ is transitive.
$(X,T)$ is called \emph{(strong) mixing} if for every pair $V, W$ of non-empty open sets in $X$, there is a compact  $K \subset T$ such that $t(V) \cap W$ is non-empty for all $t \in T \setminus K$.

 A  (semi)cascade $(X,f)$ is said to be \emph{weakly mixing} if the product system $(X \times X, f \times f)$ is transitive. $(X,f)$ is called \emph{mixing} if for every pair $V, W$ of non-empty open sets in $X$, there is a  $N>0$ such that $f^n(V) \cap W$ is non-empty for all $n \geq N$.

\br  If $(X,T)$ [$(X,f)$] is  mixing then it is weakly mixing.\er

\br Note that  $(X,T)$   is  mixing if and only if for every pair of non-empty, open
sets $U, V \subset X$ the set $N_T(U,V) = T \setminus K$, for some compact $K \subset T$. The cascade $(X,f)$ is mixing if and only if for every pair of non-empty, open sets $U, V \subset X$,  $N(U,V)$ is cofinite. \er

\br Note that for any non-empty, open $U_1, U_2, V_1, V_2 \subset X$, we have $N_T(U_1 \times V_1, U_2 \times V_2) = N_T(U_1, U_2) \cap N_T(V_1, V_2)$. \er

This observation enables us to say that for  a (semi)flow  $(X,T)$,

\bp \cite{ak} The
following are equivalent.
\be
\item  $(X,T)$ is weak mixing.

\item For non-empty,  open sets $U_1, V_1, U_2, V_2 \subset X$, there
exists a  $t \in T$ such  that $tU_1 \cap V_1 \neq
\emptyset$ and $tU_2 \cap V_2 \neq \emptyset$.

\item For non-empty,  open sets $U_1, V_1, U_2, V_2 \subset X$, there
exists a $t \in T$ such  that $t^{-1}U_1 \cap V_1 \neq
\emptyset$ and $t^{-1}U_2 \cap V_2 \neq \emptyset$.

\item The collection  $\{ N_T(U,V) : $ for $U, V$ non-empty, open in $X \}$ has the finite intersection property (or equivalently it generates a filter of subsets of $T$).

\ee

Further, if $T$ is Abelian then  for every $N \in \N$ the product system $(X^N, T)$ is
topologically transitive.

\ep

\medskip

Note that for Abelian  $T$  the product system $(X^N, T)$ being
topologically transitive is a well-known consequence of the \emph{Furstenberg Intersection Lemma}.

\bl [Furstenberg Intersection Lemma \cite{ak}] For a system $(X,T)$, with $T$ Abelian, assume that $N_T(U,V) \cap N_T(U,U)$ $ \neq  \emptyset$ for every pair of non-empty, open  $U,V \subset X$. Then for all non-empty, open  $U_1,V_1,U_2,V_2 \subset X$ there exist non-empty, open  $U_3,V_3 \subset X$ such that

$$ N_T(U_3,V_3) \quad \subset \quad  N_T(U_1,V_1) \cap  N_T(U_2,V_2). $$

 \el

\medskip

The concept of weakly mixing was first defined by \textbf{Furstenberg} in \cite{f}, who showed that the (semi)cascade $(X,f)$ is weakly mixing if each $ N(U,V)$  is thick, and is  mixing if each $ N(U,V)$  is cofinite; for  non-empty, open subsets $\ U, V$ of $X $.

\medskip

\bp \label{iweq} For a metric semicascade $(X,f)$, where $X$ is a compact and perfect,   the semicascade is weakly mixing if and only if there exists a sequence $(x_1, x_2, \ldots) \in X^\infty$ such that for all $n \in \N$, $(x_1, \ldots, x_n)$ is a transitive point for $(X^n,f^{(n)})$.
	
	We call such a sequence  a \emph{strongly transitive sequence}.
	
\ep

\bo The semicascade $(X,f)$ is weakly mixing if and only if $(X^n,f^{(n)})$ is transitive for every $n \in \N$ if and only if there exists $(x_1, \ldots, x_n) \in X^n$ that is a transitive point for $(X^n,f^{(n)})$.

We note that since  $(X^n \times X,f^{(n) } \times f)$ is transitive, the transitive point $(x_1, \ldots, x_n)$  for $(X^n,f^{(n)})$ gives $(x_1, \ldots, x_n, x_{n+1}) \in X^{n+1}$ which is a transitive point for $(X^{n+1},f^{(n+1)})$. Thus the result follows by induction. \eo

\br The existence of strongly transitive sequence is equivalent to the system being weakly mixing. This is similar to the concept of independent set due to Iwanik \cite{iw}. \er

For more on weakly mixing and  mixing, we refer to \cite{aan, f, join, gh, v1}.

\medskip

 The (semi)flow $(X,T)$ is \emph{ weakly mixing} if for all $U, V$ non-empty, open in $X $, $ N_T(U,V)  $  is thick. And is mixing if for all $U, V$ non-empty, open in $X $, $ N_T(U,V)  $  is co-compact (complement of a compact subset of $T$).

\medskip

Let $\cF$ be  a Furstenberg family. Then $(X,T)$ is \emph{$\cF$- weakly mixing} if $X \times X$ is $\cF$-transitive if and only if $\cF$ has the finite intersection property if and only if $\cF$ is thick.

A more detailed description on this can be seen in \cite{ak}.

\medskip

We recall the concept of $\cF-$point transitive from the previous section. In particular for cascades or semicascades $(X,f)$, we have;

\bp \cite{j} If $(X,f)$ is weakly mixing then  $(X, f)$   is  \emph{$\cF-$point transitive} with $\cF$  an IP family in $\Z_+$. \ep

We note that this can be proved similarly, in general, for any flow or semiflow.

\bp \label{ip} For Abelian $T$, if $(X,T)$ is weakly mixing then $(X, T)$   is  \emph{$\cF-$point transitive} with $\cF$  an IP family in $T$.\ep

Using Lemma \ref{iphit} the proof follows similarly as in \cite{j}  and we skip it.

\bd  Let $\cF$ be  a Furstenberg family such that $\cF \cdot \cF = \cF$. A sequence $(x_1, x_2, \ldots) \in X^\infty$ is   called \emph{$\cF-$strongly transitive sequence} if for all $n \in \N$, $(x_1, \ldots, x_n)$ is an  $\cF-$ transitive point for $(X^n,T)$.\ed

We have an analogue to Proposition \ref{iweq}:

\bp For Abelian $T$ and $X$  a compact and perfect metric space, if the system $(X,T)$ is weakly mixing then there exists a $\cF-$strongly transitive sequence $(x_1, x_2, \ldots) \in X^\infty$ with $\cF$  an IP family in $T$.
	\ep
	
	\bigskip

\end{subsection}

\begin{subsection}{Locally Eventually ONTO}

\bd A  semiflow $(X,S)$ is called \textbf{locally eventually onto} when for any non-empty, open $U \subset X$ there exists a compact $K \subset S$ such that $s(U) = X$, $\ \forall s \in S \setminus K$.

A  semicascade $(X,f)$ is \textbf{locally eventually onto} if and only if for any non-empty, open $U \subset X$ there exists $N \in \N$ such that $f^N(U) = X$. \ed

\br We note that for a homeomorphism $h$, it is impossible for a proper, non-empty open $U \subset X$ to satisfy $h(U) = X$. Thus, no flows or cascades can be locally eventually onto. \er

Recall Example \ref{tm} of the tent-map. This is an example of a locally eventually onto semicascade.

Also Example \ref{leo} is an example of a locally eventually onto semiflow.

Recall that $(X,f)$ is locally eventually onto then $N(U,x)$ is co-finite in $\N$ for all non-empty, open $U \subset X$ and all $x \in X$.

We recommend \cite{aan}, and the references therein, for more details on locally eventually onto semicascades.

\medskip

\bp \label{tleo} We have for a metric semiflow $(X,S)$,  the
following to be  equivalent.
\be
\item  The semiflow is locally eventually onto.

\item  For all $\epsilon > 0$, there exists a compact $K \subset T$ such that
$t^{-1}(x)$ is $\epsilon-$dense in $X$ for every $x \in X$ and $\ \forall t \in T \setminus K$.
\ee
\ep

We note that for a locally eventually onto semiflow $(X,S)$,  for every non-empty, open $U \subset X$, there exists a compact $H \subset T$ such that $\bigcup \limits_{h \in H} hU = X$.

\bc A locally eventually onto system is always mixing and hence weakly mixing and thus transitive. \ec

\bex Recall the example of the full $2-$shift.  It is mixing but fails to be locally eventually onto, as observed in \cite{aan}.\eex

$(X,S)$ is locally eventually onto then $N_S(U,x)$ is co-compact in $S$, i.e. with compact complement, for all non-empty, open $U \subset X$ and all $x \in X$.

\br A locally eventually onto system need not be almost open, as observed in \cite{aan}. \er

\bigskip

\end{subsection}

\begin{subsection}{Strongly Transitive \& Very Strongly Transitive}

The concept of transitivity deals with denseness of some orbit, while the concept of minimality implies that every  orbit is dense. What would result if we want every ``backward orbit'' to be dense? We recall a few basics on this from \cite{aan, nk} and discuss further details on this.

\medskip

 For the cascade or semicascade $(X,f)$, \emph{the backward(negative) orbit} of $x \in X$ is denoted as $\cO^{-}(x)$ and defined as,

$$\cO^{-}(x)  = \{y \in X : f^n(y) = x \  \textrm{for some}  \ n \in \N \}$$

 A cascade or semicascade $(X,f)$  is \emph{strongly transitive} if  $X = \bigcup \limits_{n \in \N} f^n(U)$  for
every non-empty, open $U \subset X$.

%
%
%
%
%
%
%
%
%
%

Let $(X,f)$ be any system, then $f$ is called \emph{iteratively almost open} if for every non-empty, open $U \subset X$
$f^n(U)^{\circ} \neq \emptyset$ for infinitely many $n \in \mathbb{N}$.

Recall that if $(X,f)$ is strongly transitive, then $f$ is iteratively almost open \cite{aan}. 

$(X,f)$ is called \emph{Very Strongly Transitive (VST)} if for every non-empty, open $U \subset X$ there is a
 $N \in \N$ such that $\bigcup \limits_{n=1}^N \ f^n(U) = X$.

Recall that if $(X,f)$ is very strongly transitive then  for every  non-empty, open set $U \subset X$ and every point $x \in X$,
 the set $N(U,x)$ is syndetic and for any non-empty, open $U, V \subset X$,
the set $N(U,V)$ is syndetic. And that if $f$ is an open map then  $(X,f)$ is very strongly transitive if and only if $(X,f)$ is strongly transitive \cite{aan}.

We refer to \cite{aan} for many other properties of strongly transitive and very strongly transitive semicascades.

\textbf{What can we say about this concept for some flow  $(X,T)$? }

Now $T^{-1}(x)  = \{y \in X : t(y) = x \  \textrm{for some}  \ t \in T \}$ is the same as $Tx$. And so the concept that every backward orbit is dense is same as the concept that every orbit is dense which is the same as minimality. So thinking about this concept for group actions gives nothing new. It just means minimality. This is also observed in \cite{aan} - for cascades strongly transitive is equivalent to being minimal, and is a distinct property only for semicascades.

\medskip

\textbf{So do we get anything new if we look into semigroup actions?}

Let $S$ be a semigroup and consider the semiflow $(X,S)$. Define $S^{-}(x)  = \{y \in X : s(y) = x \  \textrm{for some}  \ s \in S \}$.

\bd We call a  semiflow $(X,S)$   \emph{strongly transitive} if  $X = \bigcup \limits_{s \in S} sU$, for every non-empty, open set $U \subset X$.

We call a  semiflow $(X,S)$   \emph{very strongly transitive} if for every non-empty, open set $U \subset X$ there exists a compact $K \subset S$ for which $X = \bigcup \limits_{s \in K} sU$.\ed

\bp For a semiflow $(X,S)$  the following are equivalent:
\be
\item  $(X,S)$ is strongly  transitive.

\item For every non-empty, open set $U \subset X$ and every point $x \in X$, there
exists $s \in S$ such  that $x \in s(U)$.

\item  For every  non-empty, open set $U \subset X$ and every point $x \in X$,
 the hitting set $N_S(U,x)$ is non-empty.

\item $S^{-}(x)$ is dense in $X$ for
every $x \in X$.

\ee

This further implies that each point in $X$ is strongly non-wandering i.e. $\widehat{\Omega}(S) = X$.

Also, if $(X,S)$ is strongly  transitive, then $(X,S)$ is topologically transitive.

\ep

\bp If $(X,S)$ is very strongly  transitive then

\be

\item for every  non-empty, open set $U \subset X$ and every point $x \in X$,
 the hitting set $N_S(U,x)$ is syndetic, i.e. there exists a compact $K \subset S$ such that $S = K N_S(U,x)$.

\item for any non-empty, open $U, V \subset X$,
the set $N_S(U,V)$ is syndetic, i.e. there exists a compact $K \subset S$ such that $S = K N_S(U,V)$.
\ee

Further, if $(X,S)$ is very strongly  transitive then it is also strongly transitive.
\ep

Recall the definition of open semiflows from Definition \ref{open}.

\bp If semiflow $(X,S)$ is open then  the two concepts of strongly transitive and very strongly transitive are the same. \ep

We skip the trivial proofs here.

\br The concepts of strongly transitive and very strongly transitive are distinct for semicascades, as can be seen from the example in \cite{aan}. \er

\medskip

\br For   semiflow $(X,S)$, by Proposition \ref{mini} minimality implies very strongly transitive and hence strongly transitive when $S$ is Abelian and central. \er

\br For   semiflow $(X,S)$, by Theorem \ref{tleo} locally eventually onto implies very strongly transitive and hence strongly transitive. \er

\medskip

 We note that strongly transitive neither implies  nor is implied by minimality in case of semiflows. Consider the examples:

\bex  Consider  $X = \R \cup \{\infty\}$,  the one-point compactification of the  reals.  For every $r \in \R$ define $f_r: X \to X$ as $f_r(y) = y + r$, $\ \forall y \in X$. Then $S = \{f_r: r \in \R\}$ is an equicontinuous family in $\cC(X)$ (the space of all continuous real valued functions on $X$ with the uniform topology), which is a semigroup with the operation of composition of functions.

For the semiflow $(X,S)$, one has $\overline{Sx} = X$ for all $x \in X$ and so $(X,S)$ is minimal.

However, $S^-(\infty) = \emptyset$ and so $(X,S)$ is not strongly transitive.
\eex

\bex \label{leo} Let $X = \T$, and consider the semigroup $S = \R^+$ of all non-negative real numbers, with the operation addition. Consider the action $\pi: S \times X \to X$ given by $\pi(r,\theta) = \pi^r(\theta) = r\theta$ for all $r \in S$ and $\theta \in \T^1$. The semiflow $(X,S)$ is strongly transitive since $\overline{S^{-}(\theta)} = X$ for all $\theta \in X$. In fact, $S^{-}(\theta) = X$ or $S^{-}(\theta) = X \setminus \{0\}$. But $(X,S)$ is not minimal since $\{0\}$ is a minimal subset of $X$.

Also $\pi^0X = \{0 \cdot \theta: \theta \in X\}$ has empty interior and so $(X,S)$ is not almost open or open, yet $(X,S)$ is very strongly transitive.

In fact, $(X,S)$ is locally eventually onto and so also mixing and weak mixing.\eex

\medskip

For the Furstenberg family $\cF$, the semiflow $(X,S)$ is called \emph{$\cF$- strongly transitive} if for every
  non-empty, open $U \in X$ and $x \in X$, one has  $N_S(U, x) \in \cF$.

\br If $(X,S)$ is \emph{$\cF$- strongly transitive} then it is also \emph{$\cF$- transitive}. \er

\br Note that if $(X,S)$ is  \emph{$\cF$- strongly transitive} then $(X,S)$ is also \emph{$\cF$-strongly central} and  $\widehat{\Omega}_\cF(S)= X$. \er

\bigskip

\end{subsection}

\begin{subsection}{Strongly Product Transitive}

We note that  transitivity is not preserved by taking products. One can just consider the irrational rotation on $\T$ as an example. So it becomes a natural question as to when can transitivity be preserved under products. This leads to the concept of \emph{weak mixing} in topological dynamics.

\textbf{What can be said about strongly transitive  systems?}

\bex  \label{pnst} Consider  a modification of Example \ref{nao}, where we consider 

$$S = \{f_{p,q}(x) := px+q \ (\mod 1): p \in \D, q \in \D\}$$

with $\D = \{ \frac{m}{2^n}: n \in \N$ and $0 \leq m \leq 2^n\}$ and operation as composition of functions. The monoid $ S $ acts on $ X $ and $ (X,S) $ is a semiflow.

For the semiflow $(X,S)$, we have $\overline{Sx}  = X$ for all $x \in X$ and so $(X,S)$ is minimal. Note that $(X,S)$ is neither central, nor is Abelian.

Also,  $(X,S)$ is  strongly transitive since $X = \bigcup \limits_{s \in S} sU$, for every open $U \subset X$. 

 Here $S^{-} (0,0)$ is not dense in $X \times X$, and so the product semiflow $(X \times X, S)$ fails to be be strongly transitive. 

\eex

 This leads to the concept of ``strongly product transitive''. We note that strongly product transitive systems were first defined and studied in \cite{aan}. We extend that study here.

\medskip

We recall the following from \cite{aan}:

 $(X,f)$ is called \emph{Strongly Product Transitive } if for every positive integer $k$ the product system $(X^k, f^{(k)})$ is strongly transitive.

Note that for homeomorphisms, minimality is equivalent to strongly transitive, and hence homeomorphisms can never be strongly product transitive.

We recall the properties of strongly product transitive semicascades from \cite{aan}.

\bigskip

We recall the concept of multi-transitivity. $(X, f )$ is called  \emph{multi-transitive} if for each $m \in \N$, $(X^m, f \times f^2 \times \ldots \times f^m)$ is transitive \cite{tksm}. While weakly mixing means weak disjointness of $f$ with itself, multi-transitivity means weak disjointness of $f$ with its iterates. Kwietniak and Oprocha \cite{ko} have shown that in general there is no connection between weak mixing and multi-transitivity citing examples of spacing subshifts which possess just one of these properties, though they are equivalent for minimal cascades. Also weakly mixing (semi)cascades which are very strongly transitive are multi-transitive.

$(X, f )$ is called  \emph{multi-minimal} if for each $m \in \N$, the system $(X^m, f \times f^2 \times \ldots \times f^m)$ is minimal.  Examples of multi-minimal systems are discussed in \cite{ko}. We see that the concept of multi-minimality and strongly product transitive are very distinct, since multi-minimality implies that $(X^m, f \times f^2 \times \ldots \times f^m)$ is strongly transitive whereas for strongly product transitive we require $(X^m,  f^{(m)})$ is strongly transitive.

\bd A semiflow $(X,S)$ is called \emph{Strongly Product Transitive } if for every  $k \in \N$ the product system $(X^k, S)$ is strongly transitive.\ed

\br Consider the semiflow $(X,S)$. Note that for any non-empty, open $U_1, U_2,  \subset X$ and $x,y \in X$, we have $N_S(U_1 \times  U_2, x \times y) = N_S(U_1, x) \cap N_S(U_2, y)$. 

 If $U \cap V = \emptyset$, then $N_S(U, x) \cap N_S(V, x) = \emptyset$. Thus, we cannot have something like Furstenberg Intersection lemma for these hitting sets. 
 
 We can consider the transitive semicascade on an interval $(I,f)$, with $(I, f^2)$ not transitive. Look for such examples and their properties in \cite{nk}. Then we  have open intervals $U,V \subset I$, such that $f(U) = V$ and $f(V) = U$. Thus for any $x \in U $, $N(U,x)$ will comprise of even entries and $N(V,x)$ will comprise of odd entries. \er

\bp For semiflow $(X,S)$  the following are equivalent:
\be
\item  $(X,S)$ is strongly product transitive.

\item The collection  $\{ N_S(U,x) : x \in X $ for $U$ non-empty, open in $X \}$  generates a filter of subsets of $S$.
\ee \ep

The proof is trivial and follows from definition.

\bc  Strongly product transitive semiflows are weakly mixing. \ec

\bo The proof follows since for $x \in V$, $N_S(U,x) \subset N_S(U,V)$. \eo

\medskip

\br If $(X,S)$ is locally eventually onto then it is strongly product transitive.\er

\medskip

 Let $\cF$ be  a Furstenberg family.  A semiflow $(X,S)$ is called \emph{$\cF-$Strongly Product Transitive } if for every  $k \in \N$ the product system $(X^k, S)$ is $\cF-$ strongly transitive.

\br Note that for any non-empty, open $U_1, U_2  \subset X$ and $x,y \in X$, since $N_S(U_1 \times  U_2, x \times y) = N_S(U_1, x) \cap N_S(U_2, y)$,
$(X,S)$ is \emph{$\cF$- strongly product transitive} then $(X,S)$ is \emph{$\cF$- weakly mixing}  if and only if $\cF$ has the finite intersection property.\er

\bp Let $\cF$ be  a Furstenberg family such that $\cF \cdot \cF = \cF$. If $S$ is Abelian, $(X,S)$ is  \emph{$\cF$-weakly mixing} and \emph{$\cF$- strongly central} then $(X,S)$ is \emph{$\cF$- strongly product transitive}. \ep

\bo     Let $x,y \in X$ and open $U,V \subset X$ be such that $x \in U$ and $y \in V$.  For any non-empty, open $A,B \subset X$, we see that $N_S(A \times B, U \times V) \in \cF \cdot \cF = \cF$ and $N_S(U \times V, x \times y) = N_S(U,x) \cap N_S(V,y) \in \cF$.

Now $N_S(A \times B, x \times y) = N_S(A \times B, U \times V) \cdot N_S(U \times V, x \times y)$ and so is in $\cF$.

Without loss of generality, this can be extended to any finite product.

\eo

\end{subsection}

\end{section}

\begin{section}{Saransh}

We have these results for semiflows similar to the results for semicascades in \cite{aan}. We skip their trivial proofs.

\begin{theorem} Let $\phi :(X,T) \to (Y,T)$ be a factor map of dynamical systems.

 If $(X,T)$ is strongly transitive, very strongly transitive, strongly product transitive or locally eventually onto then $(Y,T)$ satisfies the corresponding property.

\end{theorem}

\begin{theorem} We have the following:

 (a) If $(X \times Y, T \times S)$ is
strongly transitive, strongly product transitive or locally eventually onto
then both $(X,T)$ and $(Y,S)$ satisfy the corresponding property.

(b) Assume $(Y,S)$ is mixing. If $(X,T)$ topologically transitive, weak mixing, or mixing then $(X \times Y, T \times S)$ satisfy the corresponding property.

(c) Assume $(Y,S)$ is locally eventually onto.  If $(X,T)$  is strongly transitive,  very strongly transitive, strongly product transitive or locally eventually onto then $(X \times Y, T \times S)$
satisfies the corresponding property. \end{theorem}

\medskip

All these properties defined above are related.

\vskip .3cm

\centerline{\scriptsize{Locally Eventually Onto $\begin{array}{l}
                                \Longrightarrow \textrm{Mixing} \\
                                \Longrightarrow \textrm{Strongly Product Transitive}
                                  \end{array}$
  $\Longrightarrow$ Weak Mixing $\Longrightarrow$ Transitive }}

  \vskip .5cm

\centerline{\scriptsize{Minimal $\Longrightarrow$ Transitive $\Longleftarrow$ Strongly Transitive $\Longleftarrow$ Very Strongly Transitive   $\Longleftarrow$ Locally Eventually Onto}}

\vskip .5cm

Moreover,

\centerline{\scriptsize{Minimal $\stackrel{\textrm{Abelian + Central}}{\Longrightarrow}$ \   \  Very Strongly Transitive $\Longrightarrow$ Strongly Transitive}}

\vskip .5cm

The reverse implications do  not hold here.

 \medskip

We try to look into those transitivies which are not compared in the above relationship implications.

\medskip

\textsc{1. Locally Eventually Onto \textit{vs} Minimality:} Consider the one-sided shift system. This semicascade is locally eventually onto but not minimal. Again Example \ref{nao} gives a Minimal semiflow which is not Locally Eventually Onto.

\medskip

\textsc{2. Minimality \textit{vs} Strongly Product Transitive:} Again the one-sided shift system gives a Strongly Product Transitive semicascade which is not minimal. Whereas Example \ref{pnst} is a Minimal semiflow which is not Strongly Product Transitive.

In fact, more can be said here. Since minimality is equivalent to strongly transitive for flows and cascades, and product of minimal systems is not minimal - these systems can never be strongly product transitive. For semicascades, minimality is almost one-to-one and guarantees the existence of an invariant subsytem which is a cascade. And so such systems cannot be strongly product transitive.

\textit{Can we have a minimal semiflow that is strongly product transitive?}

The answer is affirmative which can be seen in Example \ref{dai} below.

\medskip

\textsc{3. Strongly Product Transitive \textit{vs} Mixing:} The one-sided shift system gives a Mixing cascade which is not Strongly Product Transitive.

\textit{What can be said conversely?}

Strongly Product Transitive systems are always weakly mixing. And so if there exists a counterexample here, we are looking out for the possibility of  a weakly mixing system which is not mixing. If there exists no such counterexample, then we have a possibility that strongly product transitivity will always imply mixing.

\smallskip

\textbf{This is one of the unsolved problems mentioned in \cite{aan}. }

\smallskip

For the  case of cascades, we rule out any homeomorphic example. But all known examples  of weakly mixing but not mixing cascades are (minimal)homeomorphisms and we are not sure if a non homeomorphic (non minimal) case would exist.  \textit{It remains still open whether a semicascade example can be found in this case.}

Again, we cannot have such an example for flows. \textit{What happens in case of semiflows?} We have an affirmative answer here in form of Example \ref{dai} below.

\medskip

\medskip

\bex \label{dai} Define $\tau_q: \T \to \T$ as $\tau_q(\theta) = q \theta \ (\mod 1)$, $q \in \Q$ with $q \geq 2$ and $\mu: \T \to \T$ as $\mu(\theta) = \theta +  \alpha \ (\mod 1)$ where $\alpha \in (0,2^{-100})$ is irrational. Then $ \{\tau_q: q \in \Q, \ q \geq 2\} \cup \{\mu^m: m \in \Z\}$ generates a semigroup $S$ of continuous functions on $\T$ endowed with the  uniform topology and with the operation of composition of functions. We see that $S$ is not Abelian.

\medskip

Consider the semiflow $(\T, S)$. Note that  $(\T, S)$ is central.

\smallskip

Also $S \theta$ is dense in $\T$ for all $\theta \in \T$. Thus $(\T, S)$ is minimal.

\smallskip

Consider $S^-(\theta)$. We note that $S^-(\theta)$ is dense in $\T$  $ \ \forall \theta \in \T$ and hence $(\T, S)$ is strongly transitive.

\smallskip

Again for given non-empty, open $U \subset X$, there is a finite  $K \subset \{\mu^m: m \in \Z\}$, for which $\T = \bigcup \limits_{s \in K} \ sU$. Thus $(\T, S)$ is very strongly transitive.

\smallskip

Furthermore, we observe that  for every non-empty, open $U \subset \T $  there exists a $L \in \Q, L \geq 2$ such that $\tau_l(U) = \T$ for all $l \geq L$. Hence for all $n \in \N$ and   non-empty, open $U_1,  \ldots, U_n \subset \T$  there exists  $M_n \in \Q$  for which $\tau_{M_n}(U_1 \times \ldots \times U_n) = \T^n$.  This validates that $(\T, S)$ is strongly product transitive.

\smallskip

 We take open $U,V \subset \T$ with diameter less than $\alpha$ such that $\frac{\pi}{2} \in U$ and $\frac{3\pi}{2} \in V$. Then there exists $P \in \Q$ such that $\tau_p(U) \cap V = \emptyset$, for $p \in [2,P] \cap \Q$.
 
 Suppose $(\T,S)$ is mixing. Then for the non-empty, open $U,V \subset \T$ taken above, there exists a compact $K \subset S$ for which $\lambda(U) \cap V \neq \emptyset$ for all $\lambda \in S \setminus K$. Clearly, $\{ \tau_p: p \in [2,P] \cap \Q\} \subset K$. Let $i \in [2,P] \cap (\R \setminus \Q)$ be such that $p_n \to i$ for a sequence $\{p_n\}$ with $p_n \in [2,P] \cap \Q$. Consider $\tau_i: \T \to \T$ as $\tau_i(\theta) = i \theta \ (\mod 1)$. Then $\tau_{p_n} \to \tau_i$, but $\tau_i \notin K$ contradicting the compactness of $K$.

Hence $(\T, S)$ is not mixing.

\smallskip

\eex

\end{section}

\end{document}